\documentclass[a4paper, 12pt]{article}
\usepackage{mathrsfs}
\usepackage{amsfonts}

 \usepackage{mathrsfs,amsfonts,amsmath,amssymb}
 \usepackage{color}

\makeatletter

\@addtoreset{equation}{section}
\makeatother
\newtheorem{Theorem}{Theorem}[section]
\newtheorem{Remark}[Theorem]{Remark}

\begin{document}

\title{\textbf{The explicit solution and precise distribution of CKLS model under Girsanov transform\footnote{Supported by
China Scholarship Council, National Natural Science Foundation of
China (NSFC11026142) and Beijing Higher Education Young Elite Teacher Project (YETP0516).}}}
\author{Yunjiao Hu\ Guangqiang Lan\footnote{Corresponding author: langq@mail.buct.edu.cn. }\quad and Chong Zhang
\\ \small School of Science, Beijing University of Chemical Technology, Beijing 100029, China}

\date{}

\maketitle

\begin{abstract}
We study the relation between CKLS model and CIR model. We prove that under a suitable transformation, any CKLS model
of order $\frac{1}{2}<\gamma<1$ or $\gamma> 1$ corresponds to a CIR model under a new probability
space. Moreover, we get the explicit solution and the precise distribution of the CKLS model at any time $t$ under the
new probability measure. We also give the moment estimation of CKLS model.
\end{abstract}




\section{Introduction}

Suppose that $(\Omega,\mathscr{F},P;(\mathscr{F}_t)_{t\geq 0})$ is a right continuous filtered probability space.
In 1992, Chan, Karolyi, Longstaff and Sanders \cite{CKLS} suggested modelling the behaviour of the instantaneous
interest rate by the following stochastic differential equation
\begin{equation}\label{CKLS}\aligned dr_t&=(a-br_t)dt+\sigma r_t^\gamma dB_t, \endaligned
\end{equation}

\noindent where the initial interest rate $r_0>0,$ $a>0$, $b\in \mathbb{R}$, $\sigma>0$ and $\gamma\ge\frac{1}{2}$.
$B_t$ denotes an $1$-dimensional standard $\mathscr{F}_t$-Brownian motion.

This so called CKLS model contains many important models in finance. These models can be obtained from (\ref{CKLS})
by simply placing the appropriate restrictions on the four parameters $a, b, \sigma$ and $\gamma$. For example, if $\gamma=0$, it
becomes to a Vasicek model; $\gamma=1,a=0$ a geometric Brownian motion; $\gamma=\frac{1}{2}$ a CIR model, etc..

It is well known that there is a unique non negative solution of the equation (\ref{CKLS}).
By Mao \cite{Mao} Chapter 9, we also know that the solution $r_t>0$ for all $t\ge 0$ almost surely.
The motivation is the following:

In Lamberton and Lapeyre, \cite{LL} Proposition 6.2.5, the authors presented the Laplace transform of $r_t$ and
$\int_0^tr_sds$ in the case of CIR model (that is, $\gamma=\frac{1}{2}$). If we can transform CKLS model
into a CIR model, then we can use the known results of CIR model to the CKLS model.

Assume that $r_t$ is the unique positive solution of SDE (\ref{CKLS}).

Our main result is the following

\begin{Theorem}\label{dingli} If $\gamma>1$, or $\gamma\in(\frac{1}{2},1)$ and $\frac{\gamma}{\sigma}\ge1, b>0$,
then there exist a new probability space $(\Omega,\mathscr{F},Q;(\mathscr{F}_t)_{t\geq 0})$, a $Q$ Brownian motion $\bar{B}_t$
and a transform $f(x)=\frac{C^2}{4(1-\gamma)^2}x^{2(1-\gamma)}$ such that $Y_t=f(r_t)$ satisfies the following CIR model:
\begin{equation}\label{CIR}\aligned dY_t
=(\frac{\sigma^2 C^2}{4}+2b(1-\gamma)Y_t)dt+\sigma C\sqrt{Y_t}d\bar{B}_t.
\endaligned\end{equation}

\end{Theorem}

\begin{Remark}
When $\gamma\in(\frac{1}{2},1)$ and $\frac{\gamma}{\sigma}\ge1, b>0$, equation (\ref{CIR}) is not the usual mean
reverting square root process since the coefficient of $Y_t$ is $2b(1-\gamma)>0.$
\end{Remark}

In general, we can not solve equation (\ref{CKLS}) explicitly except for some special cases, but by a measure
transform, we can get the explicit solution of our CKLS model under the new probability space $(\Omega,\mathscr{F},Q)$.

\begin{Theorem}\label{dingli2}
Assume that $\gamma>1$, or $\gamma\in(\frac{1}{2},1)$ and $\frac{\gamma}{\sigma}\ge1, b>0$.
Under the probability measure $Q$, the explicit solution of equation (\ref{CKLS}) is
$$r_t=\Big|r_0^{1-\gamma}e^{-(\gamma-1)bt}+\sigma |\gamma-1|\int_0^te^{-b(\gamma-1)(t-s)}d\bar{B}_s\Big|^\frac{1}{1-\gamma}.$$
\end{Theorem}

We can also get the distribution of $Y_t$.

\begin{Theorem}\label{dingli3}
Assume that $\gamma>1$, or $\gamma\in(\frac{1}{2},1)$ and $\frac{\gamma}{\sigma}\ge1, b>0$.
Under the probability measure $Q$, the precise density function of the solution $r_t$ of equation (\ref{CKLS}) is
$$g_t(x)=g_{\delta,\zeta}(f(x)/L)|f'(x)|/L,\ \forall x>0,$$
where $f$ is the same as Theorem \ref{dingli}, $g_{\delta,\zeta}$ is the density of the non-central chi-square law with $\delta$
degrees of freedom and parameter $\zeta$,
$$\delta=C^2,\ L=\frac{\sigma^2C^2}{8b(\gamma-1)}(1-e^{2b(1-\gamma)}),\ \zeta=\frac{8Y_0b(\gamma-1)}{\sigma^2C^2(e^{2b(\gamma-1)}-1)}.$$

\end{Theorem}

Finally, let us consider the moment estimation of $r_t$. We have
\begin{Theorem}\label{dingli4}
Assume that $1<\gamma\le\frac{3}{2}$, or $\frac{1}{2}\le\gamma<1\ \textrm{and}\ (2\gamma+1)\sigma^2\le 2a.$
Then $$\mathbb{E}\int_0^tr_s^{-2\gamma}ds\vee\mathbb{E}\int_0^tr_s^{2(\gamma-1)}ds<\infty,\ \forall t>0,$$
where $a\vee b=\max(a,b).$\end{Theorem}

The rest of the paper is organized as follows. Firstly, we study the relation between
CKLS model and CIR model in Section 2, we derive the transform under which an arbitrary
CKLS model can be transformed to be a CIR model formally. In Section 3 we will prove that the
condition of Girsanov transform is satisfied. After this, we prove Theorem \ref{dingli2}
and Theorem \ref{dingli3} in Section 4. Finally we give the moment estimation of CKLS model.

\section{The relation between CKLS model and CIR model}

Suppose $f:\mathbb{R}_+\rightarrow\mathbb{R}_+$ is a differentiable function such that
$$x^\gamma f'(x)=C\sqrt{f(x)},\ C>0.$$

Then we have
$$f(x)=\frac{C^2}{4(1-\gamma)^2}(x^{1-\gamma}+C')^2.$$
Take $C'=0$ for simplicity. Then
$$f(x)=\frac{C^2}{4(1-\gamma)^2}x^{2(1-\gamma)}.$$

We also have
$$f'(x)=\frac{C^2}{2(1-\gamma)}x^{1-2\gamma},\quad f{''}(x)=\frac{C^2(1-2\gamma)}{2(1-\gamma)}x^{-2\gamma}.$$

It's clear that $f(x)$ is strictly monotone on $(0,\infty)$, so we have
$$f^{-1}(x)=\Big|\frac{2(\gamma-1)}{C}\Big|^\frac{1}{1-\gamma}x^\frac{1}{2(1-\gamma)}.$$

By It\^o's formula, we have
$$\aligned df(r_t)&=[f'(r_t)(a-br_t)+\frac{\sigma^2}{2}f{''}(r_t)r_t^{2\gamma}]dt+\sigma f{'}(r_t)r_t^{\gamma}dB_t\\&
=[\frac{aC^2}{2(1-\gamma)}r_t^{1-2\gamma}+\frac{bC^2}{2(1-\gamma)}r_t^{2-2\gamma}+\frac{\sigma^2}{2}
\frac{C^2(1-2\gamma)}{2(1-\gamma)}]dt+\sigma C\sqrt{f(r_t)}dB_t.
\endaligned$$

Denote $Y_t=f(r_t).$ Then $r_t=f^{-1}(Y_t)=\Big|\frac{2(\gamma-1)}{C}\Big|^\frac{1}{1-\gamma}Y_t^\frac{1}{2(1-\gamma)}.$ Therefore,
\begin{equation}\label{ito1}\aligned dY_t&=[\frac{\sigma^2 C^2(1-2\gamma)}{4(1-\gamma)}+\frac{bC^2}{2(1-\gamma)}(\frac{2(\gamma-1)}{C}\sqrt{Y_t})^2+
\frac{aC^2}{2(1-\gamma)}r_t^{1-2\gamma}]dt+\sigma C\sqrt{Y_t}dB_t.
\endaligned\end{equation}
Define

\begin{equation}\label{qt}\aligned q(r_t)&=(\frac{a}{\sigma}r_t^{-\gamma}-\frac{\gamma\sigma}{2}r_t^{\gamma-1})\textrm{sgn}(\gamma-1),\endaligned\end{equation}
$$\bar{B}_t:=B_t-\int_0^tq(r_s)ds$$
and
$$R_t=\exp\{\int_0^tq(r_s)dB_s-\frac{1}{2}\int_0^tq(r_s)^2ds\}.$$

It's clear that $R_t$ is an $\mathscr{F}_t$ local martingale with respect to probability measure $P$. If $\{R_t\}$ is a
real martingale with respect to $P$, then by Girsanov transform, there exists a probability measure $Q$ on $\Omega$ such that
$\frac{dQ}{dP}|_{\mathscr{F}_t}=R_t$ and $\bar{B}_t$ is an $\mathscr{F}_t$ Brownian motion with respect to $Q$.
Then $Y_t$ satisfies the following CIR model
$$\aligned dY_t
=(\frac{\sigma^2 C^2}{4}+2b(1-\gamma)Y_t)dt+\sigma C\sqrt{Y_t}d\bar{B}_t.
\endaligned$$

So the key point is to prove that $\{R_t\}$ is a true martingale.

\section{$R_t$ is a real martingale with respect to probability $P$}\vskip0.2in

To prove that $R_t$ is a real martingale, let us consider the following auxiliary equation
\begin{equation}\label{CKLS1}d\tilde{r}_t=(a-b\tilde{r}_t+ q(\tilde{r}_t)\sigma\tilde{r}_t^\gamma)dt+\sigma \tilde{r}_t^\gamma d B'_t,\ \tilde{r}_0=r_0,\end{equation}
where $B'_t$ is Brownian motion, $q$ are defined as in Section 2.
Since we have known that, with probability one, the solution of equation (\ref{CKLS})
will never leave the state space $(0,\infty)$, by Mijatovi\'c and Urusov \cite{MU} Corollary 2.2,
we only need to prove that the solution of equation (\ref{CKLS1})
will also never leave the state space $(0,\infty).$

\subsection{In case $\gamma>1$}

If $\gamma>1,$ then equation (\ref{CKLS1}) could be transformed to
\begin{equation}\label{CKLS2}\aligned d\tilde{r}_t&=(2a-b\tilde{r}_t-\frac{\gamma\sigma^2}{2}\tilde{r}_t^{2\gamma-1})dt+\sigma \tilde{r}_t^\gamma dB'_t.\endaligned\end{equation}

For this equation, since $\mu(x)=2a-bx-\frac{\gamma\sigma^2}{2}x^{2\gamma-1}$ satisfies
$$(x-y)(\mu(x)-\mu(y))\le 0\le K(x-y)^2,\forall x,y$$
and $x^\gamma$ is local Lipschitz continuous,
by Gy\"ongy and Krylov \cite{GK} Corollary 2.7, there exists a unique solution solution of equation (\ref{CKLS1}). We now prove that
the solution $\tilde{r}_t>0, \forall t\ge0$ almost surely.  Define
$$\tau_k:=\inf\{t>0,\ \tilde{r}_t \notin(\frac{1}{k},k)\},\ k\ge k_0$$
where $k_0$ is large enough such that $r_0\in (\frac{1}{k_0},k_0).$ We only need to prove $\tau_\infty:=\lim_{k\rightarrow\infty}\tau_k=\infty.$
If not, we can choose $T>0,\ 0<\varepsilon<1$ such that $P(\tau_\infty\le T)\ge\varepsilon$. Then there exists $k_1\ge k_0$ such that
$$P(\tau_k\le T)\ge\varepsilon\quad \forall k\ge k_1.$$

Define
$$V(x)=\sqrt{x}-1-\frac{1}{2}\log x,\ x>0.$$

Then if $\tilde{r}_t>0,$ It\^o's formula yields that
$$\aligned dV(\tilde{r}_t)&=\frac{1}{2}(\tilde{r}_t^{-\frac{1}{2}}-\tilde{r}_t^{-1})
(2a-b\tilde{r}_t-\frac{\gamma\sigma^2}{2}\tilde{r}_t^{2\gamma-1})dt\\&
\quad+\frac{\sigma^2}{4}(-\frac{1}{2}\tilde{r}_t^{-\frac{3}{2}}+\tilde{r}_t^{-2})\tilde{r}_t^{2\gamma}dt
+\frac{\sigma}{2}(\tilde{r}_t^{-\frac{1}{2}}-\tilde{r}_t^{-1}) \tilde{r}_t^\gamma dB'_t\\&
=F(\tilde{r}_t)dt+\frac{\sigma}{2}(\tilde{r}_t^{-\frac{1}{2}}-\tilde{r}_t^{-1}) \tilde{r}_t^\gamma dB'_t,\endaligned$$
where
$$F(x)=\frac{1}{2}(x^{-\frac{1}{2}}-x^{-1})
(2a-bx-\frac{\gamma\sigma^2}{2}x^{2\gamma-1})
+\frac{\sigma^2}{4}(-\frac{1}{2}x^{-\frac{3}{2}}+x^{-2})x^{2\gamma}.$$

Since the coefficients of the highest and lowest order of $F(x)$ are both negative
(the coefficients are $-\frac{\sigma^2(2\gamma+1)}{8}, -a$, respectively), then
$F(x)$ is bounded, say by $K$, on $x\in (0,\infty).$ Thus,
$$dV(\tilde{r}_t)\le Kdt+\frac{\sigma}{2}(\tilde{r}_t^{-\frac{1}{2}}-\tilde{r}_t^{-1}) \tilde{r}_t^\gamma dB'_t$$
as long as $\tilde{r}_t\in(0,\infty).$ So we have
$$\mathbb{E}(V(\tilde{r}_{T\wedge \tau_k}))\le V(r_0)+ KT.$$

Set $\Omega_k=\{\tau_k\le T\}$, then $P(\Omega_k)\ge\varepsilon$ for $k\ge k_1$. Since
$$\tilde{r}(\tau_k,\omega)=k\ \textrm{or}\ \frac{1}{k},$$
then
$$V(\tilde{r}(\tau_k,\omega))\ge(\sqrt{k}-1-\frac{1}{2}\log k)\wedge(\frac{1}{2}\log k+\sqrt{\frac{1}{k}}-1).$$

Therefore,
$$ \aligned V(r_0)+ KT&\ge\mathbb{E}(V(\tilde{r}_{T\wedge \tau_k}))\ge\mathbb{E}(1_{\Omega_k}V(\tilde{r}(\tau_k,\omega)))\\&
\ge\varepsilon\big[(\sqrt{k}-1-\frac{1}{2}\log k)\wedge(\frac{1}{2}\log k+\sqrt{\frac{1}{k}}-1)\big].\endaligned$$

Letting $k\rightarrow\infty$ leads to a contradiction
$$\infty>V(r_0)+ KT\ge\infty.$$

So $\tau_\infty=\infty.$ Therefore, $R_t$ is a martingale with respect to $P.$

\subsection{In case $\frac{1}{2}\le\gamma<1$}

In this case, equation (\ref{CKLS1}) becomes to
\begin{equation}\label{CKLS3}\aligned d\tilde{r}_t&=(\frac{\gamma\sigma}{2}
\tilde{r}_t^{2\gamma-1}-b\tilde{r}_t)dt+\sigma \tilde{r}_t^\gamma dB{'}_t.\endaligned\end{equation}

Note that the method we used in Subsection 3.1 to prove the positivity of the solution
can not be used in the present case.
Since
$$\sigma x^\gamma>0, \ \forall x>0\quad \textrm{and}\ \frac{1}{\sigma^2x^{2\gamma}},
\frac{\frac{\gamma\sigma}{2}x^{2\gamma-1}-bx}{\sigma^2x^{2\gamma}}\in L_{loc}^1(0,\infty),$$
where $L_{loc}^1(0,\infty)$ denotes the class of locally integrable functions, i.e.
the functions $(0,\infty)\rightarrow\mathbb{R}$ that are integrable on compact subsets
of $(0,\infty).$ By \cite{ES,ES1} or \cite{KS} Chapter 5, Theorem 5.15, there exists
a unique in law weak solution that possibly exits its state space $(0,\infty).$

We now prove that the solution $\tilde{r}_t>0, \forall t\ge0$ almost surely by using Feller's test for explosions.

Define the scale function
\begin{equation}p(x)=\int_1^x\exp\Big\{-2\int_1^y\frac{\frac{\gamma\sigma}{2}z^{2\gamma-1}
-bz}{\sigma^2z^{2\gamma}}dz\Big\}dy.\end{equation}

Compute
\begin{equation}p(x)=\exp\{-\frac{b}{\sigma^2(1-\gamma)}\}\int_1^xy^{-\frac{\gamma}{\sigma}}
\exp\Big\{\frac{b}{\sigma^2(1-\gamma)}y^{2(1-\gamma)}\Big\}dy.\end{equation}

In the case when $\frac{\gamma}{\sigma}\ge1$ and $b>0$, it is ease to see that
$$\lim_{x\downarrow0}p(x)=-\infty\quad \textrm{and}\quad \lim_{x\uparrow\infty}p(x)=\infty.$$

By \cite{KS} Proposition 5.22,
$$P(\tau_\infty=\infty)=1,$$
where $\tau_\infty$ is defined as Subsection 3.1. That is, the solution of equation will never leave
the state space $(0,\infty)$. Therefore, $R_t$ is a true martingale. The proof of Theorem \ref{dingli} is complete.

\section{Explicit solution and precise distribution of $r_t$ under probability $Q$}

\textbf{Proof of Theorem \ref{dingli2}}
By Theorem \ref{dingli}, we know that $Y_t=f(r_t)$ satisfies the CIR model. By It\^o's formula
$$\aligned d\sqrt{Y_t}&=\frac{1}{2\sqrt{Y_t}}[(\frac{\sigma^2 C^2}{4}+2b(1-\gamma)Y_t)dt
+\sigma C\sqrt{Y_t}d\bar{B}_t]-\frac{1}{8Y_t^\frac{3}{2}}\sigma^2C^2Y_tdt\\&
=b(1-\gamma)\sqrt{Y_t}dt+\frac{\sigma C}{2}d\bar{B}_t.\endaligned$$

So $\sqrt{Y_t}$ is an Ornstein-Uhlenbeck process, whose solution is
$$\sqrt{Y_t}=\sqrt{Y_0}e^{b(1-\gamma)t}+\frac{\sigma C}{2}\int_0^te^{b(1-\gamma)(t-s)}d\bar{B}_s.$$
Here $\sqrt{x}$ should be understood as a real number $a$ such that $a^2=x$, otherwise we will get a
contradiction since the left hand side of the equation is non negative while the right hand side will be negative with
positive probability.

Thus
$$r_t=\Big|\frac{2(\gamma-1)}{C}\Big|^\frac{1}{1-\gamma}Y_t^\frac{1}{2(1-\gamma)}=\Big|r_0^{1-\gamma}e^{-(\gamma-1)bt}+\sigma |\gamma-1|\int_0^te^{-b(\gamma-1)(t-s)}d\bar{B}_s\Big|^\frac{1}{1-\gamma}.$$

We complete the proof.

\textbf{Proof of Theorem \ref{dingli3}}
Define $L=\frac{\sigma^2C^2}{8b(\gamma-1)}(1-e^{2b(1-\gamma)})$. By \cite{LL} Proposition 6.2.5 with $\mu=0$, we know that $Y_t/L$
satisfies the non-central chi-square law with $\delta$ degrees of freedom and parameter $\zeta$, where
$$\delta=C^2, \zeta=\frac{8Y_0b(\gamma-1)}{\sigma^2C^2(e^{2b(\gamma-1)}-1)}.$$

Thus the density of $Y_t/L$ is given by
$$g_{\delta,\zeta}(x)=\frac{e^{-\zeta/2}}{2\zeta{\delta/4-1/2}}e^{-x/2}x^{\delta/4-1/2}
(\frac{\sqrt{x\zeta}}{2})^{\delta/2-1}\sum_{n=0}^\infty\frac{(\sqrt{x\zeta}/2)^{2n}}{n!\Gamma(\delta/2+n)},\quad x>0,$$
where $\Gamma(x)=\int_0^\infty t^{x-1}e^{-t}dt.$

Therefore
$$g_t(x)=g_{\delta,\zeta}(\frac{f(x)}{L})\Big|\frac{f'(x)}{L}\Big|,\quad x>0.$$

We complete the proof.

\section{Moment estimations of $r_t$}

By It\^o's formula,
\begin{equation}\label{ito}\aligned dr_t^{-2\gamma}&=(-2\gamma r_t^{-2\gamma-1}(a-br_t)+\frac{2\gamma(2\gamma+1)\sigma^2}{2}r_t^{-2})dt-2\gamma\sigma r_t^{-\gamma-1}dB_t\\&
=-2a\gamma (r_t^{-2\gamma})^{1+\frac{1}{2\gamma}}dt+ 2b\gamma r_t^{-2\gamma}dt+\gamma(2\gamma+1)\sigma^2(r_t^{-2\gamma})^{\frac{1}{\gamma}}dt+M_t.\endaligned\end{equation}

\subsection{Case I: $\frac{1}{2}\le\gamma<1$}
If $\frac{1}{2}\le\gamma<1$, since
$$(x^{-2\gamma})^{\frac{1}{\gamma}}\le (x^{-2\gamma})^{1+\frac{1}{2\gamma}}+1,\ \forall x>0,$$
we have
$$\mathbb{E}(r_t^{-2\gamma})\le r_0^{-2\gamma}+\gamma(2\gamma+1)\sigma^2t+2b\gamma\int_0^t\mathbb{E}(r_s^{-2\gamma})ds$$
in case $(2\gamma+1)\sigma^2\le 2a$.

Gronwall's lemma yields that
$$\mathbb{E}(r_t^{-2\gamma})\le \Psi(t)+2b\gamma\int_0^t\Psi(s)e^{2b\gamma(t-s)}ds,$$
where $\Psi(t)=r_0^{-2\gamma}+\gamma(2\gamma+1)\sigma^2t.$ It's clear that the right hand side of the
inequality is locally integrable with respect to $t$. Therefore,
$$\mathbb{E}\int_0^t r_s^{-2\gamma}ds<\infty$$
holds for any $t>0.$

On the other hand, by It\^o's formula,
$$\aligned dr_t^{2(\gamma-1)}&=(2(\gamma-1) r_t^{2\gamma-3}(a-br_t)+\frac{2(\gamma-1)(2\gamma-3)\sigma^2}{2}r_t^{4\gamma-4})dt+2(\gamma-1)\sigma r_t^{3\gamma-3}dB_t.\endaligned$$

Notice that $2\gamma-3\le4\gamma-4<0,$ and if $(2\gamma+1)\sigma^2\le 2a$, then
$$\frac{2(\gamma-1)(2\gamma-3)\sigma^2}{2}\le2a(1-\gamma).$$

Thus $\mathbb{E}(r_t^{2(\gamma-1)})\le r_0^{2(\gamma-1)}+(\gamma-1)(2\gamma-3)\sigma^2t+2b(1-\gamma)\int_0^t\mathbb{E}(r_s^{2(\gamma-1)})ds.$

By using Gronwall's lemma again,
$$\mathbb{E}(r_t^{2(\gamma-1)})\le \tilde{\Psi}(t)+2b(1-\gamma)\int_0^t\tilde{\Psi}(s)e^{2b(1-\gamma)(t-s)}ds,$$
where $\tilde{\Psi}(t)=r_0^{2(\gamma-1)}+(\gamma-1)(2\gamma-3)\sigma^2t.$

So we have proved Theorem \ref{dingli4} in case that
$\frac{1}{2}\le\gamma<1, (2\gamma+1)\sigma^2\le 2a.$

\subsection{Case II: $1<\gamma<\frac{3}{2}$}

Note that $0<\frac{1}{\gamma}<1$ in this case. Since
$$x^\frac{1}{\gamma}\le x+1,\ \forall x>0,$$
by (\ref{ito}) we have
$$\mathbb{E}(r_t^{-2\gamma})\le r_0^{-2\gamma}+\gamma(2\gamma+1)\sigma^2t+\gamma(2b+(2\gamma+1)\sigma^2)\int_0^t\mathbb{E}(r_s^{-\gamma})ds.$$

Gronwall's lemma yields that
$$\mathbb{E}(r_t^{-2\gamma})\le \Psi(t)+\gamma(2b+(2\gamma+1)\sigma^2)\int_0^t\Psi(s)e^{\gamma(2b+(2\gamma+1)\sigma^2)(t-s)}ds,$$
where $\Psi(t)=r_0^{-2\gamma}+\gamma(2\gamma+1)\sigma^2t.$ Then
$$\mathbb{E}\int_0^t r_s^{-2\gamma}ds<\infty$$
holds for any $t>0$.

On the other hand, by (\ref{CKLS}),
$$\mathbb{E}(r_t)=r_0+at-b\int_0^t\mathbb{E}(r_s)ds.$$
Thus $\mathbb{E}(r_t)=\frac{a}{b}+(r_0-\frac{a}{b})e^{-bt}.$
Since $1<\gamma\le\frac{3}{2}$,
$$\mathbb{E}(r_t^{2(\gamma-1)})\le 1+\mathbb{E}(r_t)=1+\frac{a}{b}+(r_0-\frac{a}{b})e^{-bt}.$$

We complete the proof of Theorem \ref{dingli4}.

From this, we know that $\mathbb{E}\int_0^tq_s^2ds<\infty$ for any $t>0$ if $\frac{1}{2}\le\gamma<1$
and $(2\gamma+1)\sigma^2\le 2a$, or $1<\gamma<\frac{3}{2}$. Therefore, we have proved that $\int_0^tq_sdB_s$ is a true martingale.

\textbf{Acknowledgement} The second author would like to thank
Professor Feng-Yu Wang for useful discussions.

\end{document}